\documentclass[12pt]{article}
\usepackage{amsmath,amsthm,amssymb,epsfig,graphicx}
\font\of=msbm10 scaled 1200
\def\R{\mbox{\of R}}
\def\C{\mbox{\of C}}
\def\Z{\mbox{\of Z}}
\def\N{\mbox{\of N}}

\def\k{\kappa}
\title{ Quadratic perturbations of quadratic codimension-four centers}
 \author{Lubomir Gavrilov \\
 \normalsize \it Institut de Math\'{e}matiques de Toulouse, UMR 5219\\
 \normalsize \it Universit\'{e} Paul Sabatier (Toulouse III)\\
   \normalsize \it 31062 Toulouse, Cedex 9 \normalsize \it   France
   \\ Iliya D. Iliev\\
\normalsize \it Institute of Mathematics, Bulgarian Academy of Sciences\\
\normalsize \it Bl. 8, 1113 Sofia, Bulgaria }

\begin{document}
\newtheorem{definition}{Definition}
\newtheorem{remark}{Remark}
\newtheorem{example}{Example}
\newtheorem{theorem}{Theorem} %[section]
\newtheorem{lemma}{Lemma}
\newtheorem{proposition}{Proposition}
\newtheorem{corollary}{Corollary}

\maketitle
 \hoffset =-1truecm \voffset =-2truecm
 \font\of=msbm10 scaled 1200

 \begin{abstract}
We study the stratum in the set of all quadratic differential systems
$\dot{x}=P_2(x,y),\;\dot{y}=Q_2(x,y)$ with a center, known as the
codimension-four case $Q_4$. It has a center and a node and a rational first
integral. The limit cycles under small quadratic perturbations in the
system are determined by the zeros of the first Poincar\'e-Pontryagin-Melnikov
integral $I$. We show that the orbits of the unperturbed system are elliptic
curves, and $I$ is a complete elliptic integral. Then using Picard-Fuchs
equations  and the Petrov's method (based on the argument principle), we set an
upper bound of eight for the number of limit cycles produced from the period
annulus around the center.
%Note that the conjectured exact upper bound should be three.
 \end{abstract}

\section{ Introduction}

The conditions for a plane polynomial quadratic vector field to have a center
are known since the beginning of the last century (Dulac (1908), Kapteyn
(1912)). In the space of all quadratic systems, the systems with a center form
a union of four irreducible affine algebraic sets  :
\begin{itemize}
    \item Hamiltonian ($Q_3^H$)
    \item reversible ($Q_3^R$)
    \item generalized Lotka-Volterra ($Q_3^{LV}$)
    \item codimension-four set ($Q_4$)
\end{itemize}
 (the subscripts indicate
the co-dimension of each algebraic subset), see \.{Z}o{\l}\c{a}dek \cite{10}
for a modern proof.
 Strictly speaking, the above
classification concerns only the centers themselves. There are quadratic
systems with two centers, one of them in $Q_3^R$ and the other in $Q_3^R\cap
Q_3^{LV}$, see \cite{5}, p. 148. In the present paper we are interested in the
maximal number of limit cycles which can emerge from the periodic orbits of a
quadratic system with a center, after a small quadratic perturbation. The
corresponding number is the cyclicity of the related open period annulus. A
quadratic center is said to be generic, if it does not belong simultaneously to
two of the above algebraic sets. The cyclicity of the period annulus of a
generic center depends on the number of the zeros of the first
Poincar\'{e}-Pontryagin-Melnikov function, while in the general case a higher-order
analysis is needed. The corresponding higher-order functions were determined in
Iliev \cite{5}.

The cyclicity of the open annuli in the generic Lotka-Volterra case
($Q_3^{LV}$) has been studied by \.{Z}o{\l}\c{a}dek \cite{10}, and in the
generic Hamiltonian case ($Q_3^H$) by Horozov-Iliev \cite{hor} and Gavrilov
\cite{gav4} (for the codimension-one generic cases from the bifurcation diagram
of $Q_3^H$ see \cite{2} and the references therein). Almost nothing is known
about the generic reversible case ($Q_3^R$) and nothing about the generic
codimension-four one ($Q_4$). In the present paper we place an upper bound on
the cyclicity of the (unique) period annulus in the generic codimension four
case ($Q_4$).
\begin{theorem}
\label{th1} The cyclicity of the open period annulus surrounding the center
of any generic codimension-four plane quadratic system is less than or equal
to eight.
\end{theorem}

The conjectural exact upper bound, as it is well known, is three \cite{10,5}.
To the end of this Introduction we briefly sketch our proof. A codimension-four
plane quadratic system with a center has a well known rational first integral
$H$. Using this one may check that the generic complexified orbit $\{(x,y)\in
\mathbb{C}^2: H(x,y)=h\}$ of the system is an affine elliptic curve, and the
corresponding  Poincar\'{e}-Pontryagin-Melnikov function $I=I(h)$ is a complete
elliptic integral on it. This remarkable fact (it seems to have gone unnoticed
by the specialists) is the starting point of the paper. Still, a direct
application of standard techniques like in \cite{8,11,hor,4,gav4} does not
work. Indeed, the Poincar\'{e}-Pontryagin-Melnikov function involves differential
of the third kind with residues algebraic in the parameter $h$ (and not only
polynomial, as in the usual cases). On the other hand, it turns out that $I$
satisfies a Picard-Fuchs equation of the form $$ M_2\circ L_2\circ L_1 (I) = 0
$$ in which $L_1 = h\frac{d}{dh} -1$ and $L_2$, $M_2$ are suitable second-order
Picard-Fuchs differential operators. The role of $L_1$ and $L_2$ in our
approach is to annihilate the residues of the differentials involved in $I$. It
is easy to see that $I$ has at most as many zeros as $G=L_1(I)$ on the interval
of consideration. We shall prove further that the solution space of $L_2$ is a
Chebyshev space (Proposition \ref{L2}). This on its turn implies that if $R$ is
a function with $k$ zeros, then each solution $G$ of the non-homogeneous
equation $L_2(G)=R$ has at most $k+2$ zeros (Proposition \ref{p2}). As
$R=L_2\circ L_1(I)= L_2(G)$ is in the kernel of the second-order Picard-Fuchs
operator $M_2$, it remains to show that $R$ has at most 6 zeros which is
achieved in a standard way by making use of the Petrov method \cite{8}, see
Proposition \ref{p4}.

\section{Remarks about the codimension-four case}
\label{s2}
  In complex coordinate $z=x+iy$, the system corresponding to the generic
  codimension-four case with a center placed at the origin becomes
$$\dot{z}=-iz+4z^2+2|z|^2+\alpha\bar{z}^2,\quad \alpha\in\C\setminus \R,
\quad |\alpha|=2.$$
It is well known that the codimension-four case has a first integral of the
form
$${\cal H}(x,y)=\frac{[\varphi(x,y)]^2}{[\psi(x,y)]^3}$$
where $\varphi$ and $\psi$ are polynomials of degree tree and two,
respectively. Taking  $\alpha=b+ic$ and denoting for short
$$Y=cx-(2+b)y,\quad \k=\frac{4}{2+b}>1,$$
one can easily find the explicit form of the first integral, namely
$${\cal H}=\frac{[8y(1+Y)-\frac23(1+\k Y^3)]^2}{[1-8y+\k Y^2]^3}.$$
Therefore, $\varphi=0$ defines a trident curve and $\psi=0$ is a parabola.
Since $\varphi(0,0)=-\frac23$, $\psi(0,0)=1$, the period annulus around the
center is placed inside the domain $\Omega$ determined by $\varphi<0<\psi$.
Taking $X^2=\psi=1-8y+\k Y^2$, $X>0$, then in $(\bar{x},\bar{y})=(X,Y)$
coordinates, the system has in $\Omega$ a first integral of the form
\begin{equation}\label{2.1}
H(\bar{x},\bar{y})=\frac{\bar{x}^{-3}}{8(2-b)}\left(\frac13\k \bar{y}^3+\k
\bar{y}^2 +(1-\bar{x}^2)\bar{y}-\bar{x}^2+\frac13\right).
\end{equation}
(Compare to \cite{5}, where an additional rotation of the very initial
coordinates takes place). Although the change of the variables
$$
(x,y)\mapsto (\bar{x},\bar{y})
$$
does not define a bi-rational projective transformation (but a double ramified
covering), it induces a bi-rational map
 $$\{(x,y): \frac{[\varphi(x,y)]^2}{[\psi(x,y)]^3} =t\}\rightarrow \{(\bar{x},\bar{y}): H(\bar{x},\bar{y})=t\}$$
 It
is seen that the level sets of the first integral $H(\bar{x},\bar{y})$  are
cubic plane curves and hence they are (generically) elliptic.

For convenience, from now on, the bars over the variables $x,y$ will be
omitted. The generating function $I(t)$, whose zeroes correspond to limit
cycles in the perturbed system, are given by the following complete elliptic
integral (cf. \cite{5}, Theorem 2 (iii))
\begin{equation}
\label{2.2}
 I(t)=\int\!\!\int_{H(x,y)<t}x^{-6}[\mu_1+\mu_2y+\mu_3y^3+
\mu_4(\kappa^2y^4-x^4)]dxdy. \end{equation}
Our main purpose in this paper is to study how many zeroes the integral $I(t)$
can have in the open interval corresponding to the period annulus around
$(1,0)$.

\section{Picard-Fuchs systems}
\label{s3}
In this section we derive several equations satisfied by the entries of
(\ref{2.2}). We also apply several reductions to express these
integrals in the form we need. In particular, we calculate explicitly
$G$, $L_2$ and $R$ we mentioned in the introduction.

We begin by introducing new variables $y=y_1-1$, $h=8(2-b)t$.
Then the equation $H(x,y)=t$ where $H$ is given by (\ref{2.1}) becomes
(we will omit the subscript in $y_1$)
\begin{equation}
\label{2.3}
H(x,y,h)=\frac{\k}{3}y^3-x^2y-hx^3-(\k-1)y+\frac23(\k-1)=0.
\end{equation}
By using (\ref{2.3}), it is not hard to verify that the period annulus around
the center at $(1,1)$ corresponds to the levels $h$ from the interval
$(-\frac23, -\frac{2}{3\sqrt{\k}})$.

Let us denote for $i,j\in\Z$
$$
I_{i,j}(h)=\int\!\!\int_{H(x,y,h)<0}x^iy^jdxdy,\quad h\in
\left(-\frac23,-\frac{2}{3\sqrt{\k}}\right).
$$
 Then multiplying (\ref{2.3})
by $x^iy^{j+1}dx$, respectively by $x^{i+1}y^jdy$ and integrating,
 we obtain the identities
\begin{equation}\label{2.5}\frac{\k}{3}(j+4)I_{i,j+3}-(j+2)I_{i+2,j+1}-h(j+1)I_{i+3,j}\end{equation}
$$\hspace{3cm}-(\k-1)(j+2)I_{i,j+1}+\frac23(\k-1)(j+1)I_{i,j}=0,$$
\begin{equation}\label{2.6}\frac{\k}{3}(i+1)I_{i,j+3}-(i+3)I_{i+2,j+1}-h(i+4)I_{i+3,j}\end{equation}
$$\hspace{3cm}-(\k-1)(i+1)I_{i,j+1}+\frac23(\k-1)(i+1)I_{i,j}=0.$$
Multiplying (\ref{2.5}) by $i+4$ and (\ref{2.6}) by $j+1$ and
subtracting, we come to
$$
\k(i+j+5)I_{i,j+3}-(i+j+5)I_{i+2,j+1}-(\k-1)(i+3j+7)I_{i,j+1}
+2(\k-1)(j+1)I_{i,j}=0.
$$
 In particular, for $i=-6$
and $j=1$ one obtains $I_{-6,2}=I_{-6,1}$. Therefore, the function
in (\ref{2.2}) takes the form
\begin{equation}\label{2.8} I(h)=\mu_1I_{-6,0}+\mu_2I_{-6,1}+\mu_3I_{-6,3}
+\mu_4(\k^2I_{-6,4}-I_{-2,0})\end{equation} with all constants $\mu_i\in\R$
independent.

Let us apply to (\ref{2.3}) and (\ref{2.8}) the change of variables $(x,y) \to
(x^{-1}, yx^{-1})$. Then (\ref{2.3}) reduces to
\begin{equation}\label{2.9}H(x,y)\equiv \frac23(\kappa-1)x^3-(\kappa-1)x^2y+\frac{\kappa}{3}y^3-y=h,\end{equation}
$I_{i,j}(h)$ becomes $-I_{-i-j-3,j}(h)$ and (\ref{2.8}) becomes
\begin{equation}\label{2.10}I(h)=\int\!\!\int_{H(x,y)<h}\left(\mu_1x^3+\mu_2x^2y+\mu_3y^3+
\mu_4\left(\frac{\kappa^2y^4-1}{x}\right)\right)dxdy.\end{equation} By
(\ref{2.9}), $H(x,y)=-H(-x,-y)$, therefore the phase portrait of the related
Hamiltonian system has a central symmetry with respect to the origin.

Next, we can use the following identities:
$$\begin{array}{l}\displaystyle
I_{1,2}=I_{2,1}=\frac{3h}{10}I_{0,0}+I_{1,0}+\frac15 I_{0,1},\\[3mm]
\displaystyle
I_{3,0}=\frac{3\k h}{10(\k-1)}I_{0,0}+I_{1,0}+\frac{\k}{5(\k-1)}I_{0,1},\\[3mm]
\displaystyle
I_{0,3}=\frac{3(\k+1)h}{10\k}I_{0,0}+\frac{\k-1}{\k}I_{1,0}+\frac{\k+6}{5\k}I_{0,1},\\[3mm]
\displaystyle
I_{-1,4}=\frac{6h}{5\k}I_{-1,1}+\frac{9}{5\k^2}I_{-1,0}+\frac{9(\k-1)}{5\k^2}I_{1,0}
+\frac{\k-1}{\k}I_{1,2}
\end{array}$$
to transform (\ref{2.10}) into
\begin{equation}
\label{2.11}
I(h)=\mu_1hI_{0,0}(h)+\mu_2I_{1,0}(h)+\mu_3I_{0,1}(h)+\mu_4[2I_{-1,0}(h)+3\k
hI_{-1,1}(h)]. \end{equation} Following the standard way \cite{hor}, one can
derive a Picard-Fuchs system for the entries in (\ref{2.11}). Its explicit form
is as follows.
$$\begin{array}{l}
\displaystyle
I_{0,0}=\frac{3h}{2}I_{0,0}'+I_{0,1}',\\[3mm]
\displaystyle
I_{1,0}=hI_{1,0}'+\frac23 I_{1,1}',\\[3mm]
\displaystyle
I_{0,1}=\frac{2}{3\k}I_{0,0}'+hI_{0,1}'+\frac{2(\k-1)}{3\k}I_{1,1}',\\[3mm]
\displaystyle
I_{1,1}=\frac{3h}{8}I_{0,0}'+\frac12 I_{1,0}'+\frac14 I_{0,1}'+\frac{3h}{4}I_{1,1}',\\[3mm]
\displaystyle
I_{-1,0}=3hI_{-1,0}'+2I_{-1,1}',\\[3mm]
\displaystyle
I_{-1,1}=\frac{\k-1}{\k}I_{1,0}'+\frac{1}{\k}I_{-1,0}'+\frac{3h}{2}I_{-1,1}'.
\end{array}$$
By using the above system, we see that
$$
hI'-I=G(h),\quad G(h)=(\mu_1h^2+\mu_3)I'_{0,0}+\mu_2I'_{1,1}
+\mu_4[-4hI'_{-1,0}+(3\k h^2-4)I'_{-1,1}].
$$
Therefore
$$
I(h)=h\int_{-\frac23}^h\xi^{-2}G(\xi)d\xi
$$
and $I(h)$ has at most as much zeroes as $G(h)$ in
$(-\frac23,-\frac{2}{3\sqrt{\k}})$.
By the same system, the integrals $I_{0,0}'$ and $I_{1,1}'$ satisfy
\begin{equation}\label{pfs}
\begin{array}{lll}
-3\k h I'_{0,0}&=&(9\k h^2-4)I''_{0,0}-4(\k-1)I''_{1,1},\\
-3\k hI'_{1,1}&=&(9\k h^2-4)(I''_{0,0}-I''_{1,1}),
\end{array}
\end{equation}
and the integrals $I'_{-1,0}$ and $I'_{-1,1}$ satisfy
$$
\begin{array}{l}
\displaystyle
I'_{-1,0}=-\frac{3h}{2}I''_{-1,0}-I''_{-1,1},\\[3mm]
\displaystyle
I'_{-1,1}=-\frac{2}{\k}I''_{-1,0}-3hI''_{-1,1}+\frac{4(\k-1)}{3\k h}I''_{1,1}.
\end{array}
$$
Hence, the integral
$J=-4hI'_{-1,0}+(3\k h^2-4)I'_{-1,1}$ satisfies the second-order equation
$$L_2(h)J=\frac43(\k-1)[h(9\k h^2-4)I'''_{1,1}+(6\k h^2+8)I''_{1,1}],$$
with
$$L_2(h)=5\k h-(9\k h^2-8)\frac{d}{dh}+h(9\k h^2-4)\frac{d^2}{dh^2}.
$$
Therefore, a similar equation $L_2(h)G=R$ (with a right-hand side $R$ depending
linearly on $\mu_i$, $I'_{0,0}$ and $I'_{1,1}$ ) also holds. To calculate $R$
explicitly, we first use (\ref{pfs}) to obtain the identities
$$\begin{array}{l}
\displaystyle I_{0,0}''=\frac{-3h(9\k h^2-4)I'_{0,0}+12(\k-1)hI'_{1,1}}
{(9h^2-4)(9\k h^2-4)},\\[4mm]
\displaystyle I_{1,1}''=\frac{-3hI'_{0,0}+3hI'_{1,1}} {9h^2-4},\\[4mm]
\displaystyle I_{0,0}'''=\frac{324\k h^4+(72\k-108)h^2-48}
{(9h^2-4)^2(9\k h^2-4)}I'_{0,0}-\frac{12(\k-1)[243\k h^4-36(\k+1)h^2-16]}
{(9h^2-4)^2(9\k h^2-4)^2}I'_{1,1},\\[4mm]
\displaystyle I_{1,1}'''=\frac{27h^2+12}
{(9h^2-4)^2}I'_{0,0}-\frac{162\k h^4+(144\k-108)h^2-48}
{(9h^2-4)^2(9\k h^2-4)}I'_{1,1}.
\end{array}$$
A direct calculation then yields
\begin{equation}\label{r}
R(h)=\frac{h[(a_0+a_1h^2+a_2h^4+a_3h^6)I'_{0,0}+(b_0+b_1h^2+b_2h^4)I'_{1,1}]}
{(9h^2-4)^2(9\k h^2-4)},
\end{equation}
with some constants $a_j$, $b_j$ depending linearly on $\mu_i$.
Below, we shall use the explicit formulas for $G$, $R$, $L_2$ just derived
in order to prove our main result.

\section{Proof of Theorem \ref{th1}}
\label{s4} The proof of Theorem \ref{th1} follows from the next four
Propositions, the first two of them being probably known. Let $V$ be a
finite-dimensional vector space of functions, real-analytic on an open interval
$(a,b)$.
\begin{definition}
We say that $V$ is a Chebyshev space, provided that each non-zero function in
$V$ has at most $\dim(V)-1$ zeros, counted with multiplicity.
\end{definition}

Let $S$ be the  solutions space of a second-order linear analytic differential
equation
\begin{equation}
\label{eq}
 x'' + a_1(t) x' +a_2(t) x = 0
\end{equation}
on an open interval $(a,b)$.
\begin{proposition}
\label{ect} The solution space $S$ of $(\ref{eq})$ is a Chebyshev space on the
interval $(a,b)$ if and only if there exists a nowhere vanishing solution
$x_0\in S$ $(x_0(t)\neq 0$, $\forall t\in (a,b))$.
\end{proposition}
\begin{remark}
{\em The question of existence of a non-vanishing solution is a recurrent
question in many papers concerning zeros of Abelian integrals, see e.g.
\cite{11,gi04,coga}. A Chebyshev space $V$ in our sense is sometimes called
an \emph{extended} Chebyshev space, and it is said to be an \emph{extended
complete} Chebyshev space, provided that it has a complete flag of extended
Chebyshev sub-spaces, see e.g. \cite{karlin}. In the case when $\dim(V)=2$
the Chebyshev space $V$ (in our sense) is an extended complete Chebyshev one
if and only if it has a nowhere vanishing function. Therefore  the notions of
Chebyshev space (in our sense) and extended complete Chebyshev space (in the
sense of \cite{karlin}), as far as applied to the solution space of
(\ref{eq}), coincide. }
\end{remark}
\begin{proposition}
\label{p2} Suppose the solution space of the homogeneous equation $(\ref{eq})$
is a Chebyshev space and let $R$ be an analytic function on $(a,b)$ having $k$
zeros (counted with multiplicity). Then every solution $x(t)$ of the
non-homogeneous equation
\begin{equation}\label{nonhomeq}
 x'' + a_1(t) x' +a_2(t) x = R(t)
\end{equation}
has at most $k+2$ zeros on $(a,b)$.
\end{proposition}

\begin{proposition}
\label{pp} The solution space $S$ associated to the differential operator
\begin{equation}\label{L2}
L_2(h)=5\k h-(9\k h^2-8)\frac{d}{dh}+h(9\k h^2-4)\frac{d^2}{dh^2},\quad \k>1
\end{equation}
is a Chebyshev system on the interval $(-\infty,-\frac{2}{3\sqrt{\k}})$.
\end{proposition}

Let $R=L_2\circ L_1 (I)$, where $I$ is the Abelian integral (\ref{2.11}) and
$L_1 = h\frac{d}{dh} -1$. A suitable for our purposes expression for $R$ is
obtained in (\ref{r}).
\begin{proposition}
\label{p4} The Abelian integral $R(h)$ has at most $6$ zeros (counted with
multiplicity) on the interval $(-\frac{2}{3},-\frac{2}{3\sqrt{\k}})$, $\k>1$.
\end{proposition}
\noindent \textbf{Proof of Theorem \ref{th1}, assuming Propositions
\ref{ect}-\ref{p4}.} The Abelian integral $L_1 (I)$, $I$ given by (\ref{2.11}),
is a solution of the non-homogeneous equation $L_2(G)=R$. According to
Proposition \ref{p4}, Proposition \ref{pp} and Proposition \ref{p2} the
integral $L_1 (I)(h)$ has at most $8$ zeros on the interval
$(-\infty,-\frac{2}{3\sqrt{\k}})$. The integral $I(h)$ has the same number of
zeros as $L_1 (I)(h)$ on the same interval. Finally, the functions
(\ref{2.11}), (\ref{2.10}), (\ref{2.8}) and (\ref{2.2}) have the same number of
zeros in the respective intervals. $\Box$

To the end of the paper we prove the above Propositions 1-4. \\

\vspace{1ex}
\noindent
\textbf{ Proof of Proposition \ref{ect}.} Let $(x_1,x_2)$ be a fundamental set
of solutions of (\ref{eq}) and consider the map
%%\begin{equation}\label{map}
$$
 p:  (a,b)\rightarrow S^1=\mathbb{P}^1\mathbb{R} : t\mapsto [x_1(t):x_2(t)] .
$$
%%\end{equation}
As the Wronskian of $x_1, x_2$ is non-vanishing, then the map $p$ is
non-degenerate ($dp(t)\neq 0$) and hence monotonous.

The solution space of (\ref{eq}) is Chebyshev if and only if the map $p$ is
injective. The solution space of (\ref{eq}) contains a nowhere vanishing
solution if and only if the map $p$ is not surjective.

As the circle is not homeomorphic to an open interval, then the monotonous
differentiable map $p$ cannot be surjective and injective at the same time. It
follows that if $p$ is injective then it is not surjective. If, on the
contrary, $p$ is not surjective, then the monotonicity of $p$ implies that the
image of $(a,b)$ under $p$ is an open subinterval of $S^1$ and $p$ is
injective.
$\Box$ \\

\vspace{1ex}
\noindent
\textbf{ Proof of Proposition \ref{p2}.}  Let $(x_1,x_2)$ be a
fundamental set of solutions of (\ref{eq}), such that
 $x_1(t)$ is a nowhere vanishing solution. The change of the variables
 $x\rightarrow x/x_1(t)$ does not change the number of the zeros of the
 solutions of (\ref{eq}), which is transformed to a linear equation with a fundamental system of
 solutions $\{1,\frac{x_2(t)}{x_1(t)}\}$. As the vector space spanned by
 $x_1,x_2$ is  Chebyshev, then
the function $\frac{x_2(t)}{x_1(t)}$
 is strictly monotonous on $(a,b)$. The change of the independent variable
 $t\rightarrow \tau = \frac{x_2(t)}{x_1(t)}$ is therefore regular and
 transforms further the above linear equation to an equation with a fundamental system of
 solutions $\{1,\tau\}$. Therefore the corresponding differential operator is a multiple of $\frac{d^2}{d\tau^2}$.
 More precisely, the regular change of variables
$$
(x,t) \mapsto (y,\tau),\quad y=\frac{x}{x_1(t)},\quad \tau =
\frac{x_2(t)}{x_1(t)}
 $$
transforms equation (\ref{nonhomeq}) to
$$
x_1 \left(\frac{d}{dt} \frac{x_2(t)} {x_1(t)}\right)^2 \frac{d^2}{d\tau^2} y =
R(t(\tau))
$$
and hence each solution of the non-homogeneous equation (\ref{nonhomeq}) has at
most $k+2$ zeros on $(a,b)$ (counted with multiplicity). $\Box$ \\

\vspace{1ex}
\noindent
{\bf Proof of Proposition \ref{pp}.}
Let $\{\delta(h): h\in (-\frac23,-\frac{2}{3\sqrt{\k}})\}$
be the continuous family of periodic orbits defined by $\{H=h\}$, with $H$
in the form (\ref{2.9}). Then $G(h)=
\int_{\delta(h)} \omega$ where $\omega$ is a linear combination of elliptic
differentials of the first and second kind. Therefore the residues of $\omega$
are solutions of $L_2$. The only residues of $\omega$ are at $(0,y)$ where $y$
is one of the roots of $\frac{\k}{3} y^3-y=h$ and they are easily computed:
$$
Res_{(0,y)} \omega = \frac{-4h+ (3 \k h^2-4)y}{\k y^2-1} .
$$
For $h< -\frac{2}{3\sqrt{\k}}$ the polynomial $\frac{\k}{3} y^3-y - h$ has one
real root which we denote by $y_0$. We shall show that the solution $ f(h)=
Res_{(0,y_0)} \omega $ of $L_2x=0$ does not vanish in the interval
$(-\infty,-\frac{2}{3\sqrt{\k}})$. Indeed, on this interval $y_0(h)$ is a
strictly increasing function and $y_0(h)<
-\sqrt{\frac{5}{\k}}=y_0(-\frac{2}{3\sqrt{\k}})$. It remains to show that $
-4h+ (3 \k h^2-4)y_0 \neq 0$. The identity
$$
-4h+ (3 \k h^2-4)y_0= -4(\frac{\k}{3} y_0^3-y_0)+ (3 \k h^2-4)y_0= \k
h(3h^2-\frac{4}{3} y_0^2)
$$
implies that $f(h)=0$ on $(-\infty,-\frac{2}{3\sqrt{\k}})$ if and only if $h=
2y_0/3$. Now
$$\frac{\k}{3} y_0^3-y_0 = h = \frac{2}{3}y_0$$
gives $y_0=\pm \sqrt{\frac{5}{\k}}$. But $y_0=y_0(h)$ is a strictly increasing
function in $(-\infty,-\frac{2}{3\sqrt{\k}}]$ and $y_0(-\frac{2}{3\sqrt{\k}})=
-\sqrt{\frac{5}{\k}}$ which is the needed contradiction. Thus the solution
space of $L$ on $(-\infty,-\frac{2}{3\sqrt{\k}})$ contains a nowhere vanishing
function and hence is a Chebyshev system.  $\Box$
% proof p2

 The above result cannot be improved, as shown by the example
$L=\frac{d^2}{dt^2}$.\\

\noindent{\bf Proof of Proposition \ref{p4}.} According to (\ref{r}), it
suffices to show that any linear combination of the form
$P_3(h^2) I'_{0,0}(h) + Q_2(h^2) I_{1,1}'(h)$ where $P_3, Q_2$ are real
polynomials of degree at most three and two, has at most 6 zeros. We note that
$I'_{0,0}, I_{1,1}'$ are complete elliptic integrals of the first and second
kind respectively, satisfying the second-order Picard-Fuchs system (\ref{pfs}).

We introduce a new variable $s\in(1,\k)$
through $h=-\frac{2}{3}\sqrt{s/\k}$ and denote by dot the
differentiation with respect to $s$. Also, denote for a convenience
$J_1(s)=I'_{0,0}(h(s))$, $J_2(s)=I'_{1,1}(h(s))$.
By (\ref{L2}) and (\ref{r}), we obtain the equation
%%\begin{equation}\label{gauss}
$$
L_2G\equiv\left[s(1-s)\frac{d^2}{ds^2}-\frac12\frac{d}{ds}
-\frac{5}{36}\right]G(s)= \frac{P_3(s)J_1(s)+Q_2(s)J_2(s)}
{(s-\k)^2(s-1)}.
$$
%%\end{equation}
Hence, we will need information about the zeroes of the right-hand side
in the interval $(1,\k)\subset(1,\infty)$. Equation (\ref{pfs}) implies
that $J(s)= (J_1(s), J_2(s))^\top$ satisfies the system of hypergeometric
type
%%\begin{equation}\label{pfs1}
$$
    J(s)= 6 \left(
    \begin{array}{cc}
      1-s & \k-1 \\
      1-s & s-1 \\
   \end{array}\right) \dot{J}(s)
$$
%%\end{equation}
or equivalently
\begin{equation}\label{pfs2}
6 (s-1)(s-\k) \dot{J}(s) = \left(
    \begin{array}{cc}
      1-s & \k-1 \\
      1-s & s-1 \\
    \end{array}\right) J(s)
\end{equation}

Let us consider for any $n\in\N$ the vector space
$$V_n= \{P_n. J_1 + Q_{n-1} .  J_2: P_n, Q_{n-1} \in \mathbb{R}[h],\quad \deg P_m,
Q_m \leq m \} .
$$
Proposition \ref{p4} follows from the following more general result.
\begin{proposition}
\label{p5} The vector space $V_n$ is Chebyshev on the interval $(1,\k)$ : each
element has at most $\dim V_n - 1= 2n$ zeros (counted with multiplicity).
\end{proposition}
\noindent{\bf Proof of Proposition \ref{p5}.} We use the Petrov method in the
complex domain ${\cal D}=\mathbb{C}\setminus~(-\infty, 1)$, see \cite{8,11}. The
characteristic exponents of (\ref{pfs2}) at $1, \k, \infty$ are equal to
$\{0,0\}$, $\{0,0\}$, $\{-\frac{1}{6},\frac{1}{6}\}$, respectively. The
function $P_n. J_1 + Q_{n-1}. J_2$ is holomorphic in a neighborhood of $s=\k$
(this value corresponds to the center of the system $dH=0$, with $H$ the
symmetric Hamiltonian given by (\ref{2.9})), and has a logarithmic singularity
in a neighborhood of $s=1$ (which corresponds to the saddle point of the
symmetric Hamiltonian system). The function $J_1$ is a complete elliptic
integral of the first kind and therefore does not vanish. Consider the function
$$
F(s)=\frac{P_n(s). J_1(s) + Q_{n-1}(s) .  J_2(s)}{J_1(s)}
$$
which is real-analytic in the complex domain $\mathcal{D}$.
We apply the argument principle to the smaller domain
$$\mathcal{D}_\varepsilon = \mathcal{D}\cap \{s:
|s-1| > \varepsilon \}\cap \{s: |s| < \frac{1}{\varepsilon}\}.$$ For this
purpose, we consider the increase (or decrease) of the argument of $F$ when $s$
makes one turn along the boundary of $\mathcal{D}_\varepsilon$ in a positive
direction. The following facts are easily deduced from the asymptotic
expansions of $J$ near the singular points of the Fuchs system (\ref{pfs2}).
\begin{enumerate}
    \item Along the boundary of the small disc $\{|s-1|=\varepsilon \}$ the
    increase of the argument of $F$ is bounded by a value close to zero.
    \item Along the boundary of the
big disc $\{|s|=\frac{1}{\varepsilon } \}$ the increase of the argument of $F$
is bounded by a value close to $2\pi \, n = 2\pi \max\{n,n-1+\frac{2}{6} \}$.
    \item Along the interval $(-\infty,1)$, the imaginary part of $F$ equals
$$
Q_{n-1}(s)\; \mbox{\rm Im}\frac{J_2(s)}{J_1(s)}= Q_{n-1}(s) \frac{\det
W(s)}{|J_1(s)|^2}
$$
where
$$
W(s)=\left(
\begin{array}{cc}
  J_1(s) & \tilde{J}_1(s)\\
  J_2(s) & \tilde{J}_2(s) \\
\end{array}\right)
$$
is a fundamental matrix of (\ref{pfs2}).
\item The determinant of the fundamental matrix $W$ is a rational function in
$s$ and in fact a non-zero constant.
\end{enumerate}
Summing up the above facts we conclude that the increase of the argument of $F$
along the boundary of $\mathcal{D}_\varepsilon$ is bounded by $2n$. Therefore
$F$, and hence $P_n(s). J_1(s) + Q_{n-1}(s) .  J_2(s)$ has at most $2n$ zeros
(counted with multiplicity) in $\mathcal{D}$, and hence in $(1,\k)$.
Proposition \ref{p5}, and hence Proposition \ref{p4} are proved. $\Box$

This also finishes the proof of Theorem  \ref{th1}.

\vspace{2ex} \noindent {\bf Acknowledgment.} Part of the paper was written
while the second author was visiting the University  of Toulouse. He thanks for
its hospitality. This research has been partially
 supported by PAI  Rila program through Grants  14749SM (France) and Rila 3/6-2006 (Bulgaria).

\vspace{2cm}

E-mail addresses:

 \vspace{1cm}

lubomir.gavrilov@math.ups-tlse.fr

iliya@math.bas.bg

\end{document}